\documentclass[number,seceqn,dvips]{arxbj}

\aid{0}
\volume{16}
\issue{2}
\pubyear{2010}
\firstpage{459}
\lastpage{470}
\doi{10.3150/09-BEJ216}

\makeatletter
\newtheorem{theorem}{Theorem}

\newtheorem{lemma}{Lemma}

\makeatother

\begin{document}
\begin{frontmatter}

\title{Relative log-concavity and a pair of triangle inequalities}
\runtitle{A pair of triangle inequalities}

\begin{aug}
\author{\fnms{Yaming} \snm{Yu}\ead[label=e1]{yamingy@uci.edu}}
\runauthor{Y. Yu}
\address{Department of Statistics, University of California,
Irvine, CA 92697-1250, USA.\\ \printead{e1}}
\end{aug}

\received{\smonth{3} \syear{2008}}
\revised{\smonth{5} \syear{2009}}

%
\begin{abstract}
The \textit{relative log-concavity} ordering $\leq_{\mathrm{lc}}$ between
probability mass functions (pmf's)
on non-negative integers is studied. Given three pmf's $f, g, h$ that
satisfy $f\leq_{\mathrm{lc}} g\leq_{\mathrm{lc}} h$,
we present a pair of (reverse) triangle inequalities: if $\sum_i
if_i=\sum_i ig_i<\infty,$ then
\[
D(f|h)\geq D(f|g)+D(g|h)
\]
and if $\sum_i ig_i=\sum_i ih_i<\infty,$ then
\[
D(h|f)\geq D(h|g)+D(g|f),
\]
where $D(\cdot|\cdot)$ denotes the Kullback--Leibler divergence. These
inequalities,
interesting in themselves, are also applied to several problems,
including maximum entropy
characterizations of Poisson and binomial distributions and the best
binomial approximation
in relative entropy. We also present parallel results for continuous
distributions and discuss
the behavior of $\leq_{\mathrm{lc}}$ under convolution.
\end{abstract}

%

\end{frontmatter}

\section{Introduction and main result}\label{sec1}

A non-negative sequence $u=\{u_i, i\geq0\}$ is \textit{log-concave} if
(a) the support of $u$ is an interval in $\mathbf{Z}_+=\{0, 1,\ldots\}$
and (b) $u_i^2\geq u_{i+1}u_{i-1}$ for all $i$ or, equivalently, $\log
(u_i)$ is concave in $\operatorname{supp}(u)$. Such sequences occur naturally in
combinatorics, probability and statistics, for example, as probability
mass functions
(pmf's) of many discrete distributions. Given two pmf's $f=\{f_0, f_1,
\ldots\}$ and $g=\{g_0, g_1,\ldots\}$ on
$\mathbf{Z}_+$, we say that $f$ is log-concave relative to $g$, written
as $f\leq_{\mathrm{lc}} g,$ if
\begin{longlist}
\item[1.]
each of $f$ and $g$ is supported on an interval on $\mathbf{Z}_+$;
\item[2.]
$\operatorname{supp}(f)\subset \operatorname{supp}(g)$;
\item[3.]
$\log(f_i/g_i)$ is concave in $\operatorname{supp}(f)$.
\end{longlist}
We have $f\leq_{\mathrm{lc}} f$ (assuming interval support) and $f\leq_{\mathrm{lc}} g,\
g\leq_{\mathrm{lc}} h\Longrightarrow f\leq_{\mathrm{lc}} h$. In other words,
$\leq_{\mathrm{lc}}$ defines a pre-order among discrete distributions with
interval supports on $\mathbf{Z}_+$. When $g$ is a geometric pmf, $f\leq
_{\mathrm{lc}} g$ simply means that $f$ is log-concave; when $g$ is a binomial
or Poisson
pmf and $f\leq_{\mathrm{lc}} g$, then $f$ is \textit{ultra log-concave} \cite{P}
(see Section~\ref{sec2}).

Whitt \cite{W} discusses this particular ordering and illustrates its
usefulness with a queueing theory example.
Yu \cite{Yu2} uses $\leq_{\mathrm{lc}}$ to derive simple conditions that imply
other stochastic orders such as the usual stochastic order, the hazard
rate order and the likelihood ratio order. Stochastic orders play an
important role in diverse areas, including reliability theory and
survival analysis (\cite{BP,DJ}); see Shaked and Shanthikumar \cite
{SS} for a book-length treatment. In this paper, we are concerned with
entropy relations between distributions under $\leq_{\mathrm{lc}}$. The
investigation is motivated by maximum
entropy characterizations of binomial and Poisson distributions (see
Section~\ref{sec2}). For a random variable $X$ on $\mathbf{Z}_+$ with pmf $f$,
the Shannon entropy is defined as
\[
H(X)=H(f)=-\sum_{i=0}^\infty f_i\log(f_i).
\]
By convention, $0\log(0)=0$. The relative entropy (Kullback and Leibler
\cite{KL}; Kullback \cite{K}; Csisz\'{a}r and Shields \cite{CS})
between pmf's $f$ and $g$ on $\mathbf{Z}_+$ is defined as
\[
D(f|g)=
\cases{\displaystyle\sum_{i=0}^\infty f_i\log(f_i/g_i), &\quad  if
$\operatorname{supp}(f)\subset \operatorname{supp}(g)$,
\cr
\infty, &\quad otherwise.
}
\]
By convention, $0\log(0/0)=0$. We state our main result.
\begin{theorem}
\label{main}
Let $f, g, h$ be pmf's on $\mathbf{Z}_+$ such that $f\leq_{\mathrm{lc}} g\leq
_{\mathrm{lc}} h$. If $f$ and $g$ have finite and equal means, then
$D(f|h)<\infty$ and
%
\begin{equation}
\label{tri}
D(f|h)\geq D(f|g)+D(g|h);
\end{equation}
if $h$ and $g$ have finite and equal means, then
%
\begin{equation}
\label{tri-}
D(h|f)\geq D(h|g)+D(g|f).
\end{equation}
\end{theorem}

Theorem~\ref{main} has an appealing geometric interpretation. (With a slight
abuse of notation, we write the mean of a pmf $g$ as $E(g)=\sum_i
ig_i$.) If $g$ and $h$ satisfy $E(g)<\infty$ and $g\leq_{\mathrm{lc}} h$, then
(\ref{tri}) gives
\[
D(g|h)=\inf_{f\in F} D(f|h),\qquad F=\{f\dvt f\leq_{\mathrm{lc}} g, E(f)=E(g)\}.
\]
That is, $g$ is the \textit{I-projection} of $h$ onto $F$. Relation (\ref
{tri-}) can be interpreted similarly. See
Csisz\'{a}r and Shields \cite{CS} for general definitions and
properties of the I-projection and the related \textit{reverse
I-projection}.

While Theorem \ref{main} is interesting in itself, it can also be used
to derive several classical and new
entropy comparison results. We therefore defer its proof to Section~\ref{sec3},
after considering these applications. We conclude in Section~\ref{sec4} with
extensions to continuous distributions. Throughout, we also discuss the
behavior of $\leq_{\mathrm{lc}}$ under convolution, as this becomes relevant in
a few places.

\section[Some implications of Theorem 1]{Some implications of Theorem \protect\ref{main}}\label{sec2}

Theorem \ref{main} is used to unify and generalize classical results on
maximum entropy characterizations of Poisson and
binomial distributions in Section~\ref{sec2.1} and to determine the best
binomial approximation to a sum of independent Bernoulli
random variables (in relative entropy) in Section~\ref{sec2.2}.
Section~\ref{sec2.3}
contains analogous results for the negative binomial. Theorem~\ref{main}
also implies monotonicity (in terms of relative entropy) in
certain Poisson limit theorems.

\subsection{Maximum entropy properties of binomial and Poisson
distributions}\label{sec2.1}

Throughout this subsection (and in Section~\ref{sec2.2}), let $X_1,\ldots, X_n$
be independent Bernoulli random variables with
$\Pr(X_i=1)=1-\Pr(X_i=0)=p_i, 0<p_i<1$. Define $S=\sum_{i=1}^n X_i$
and $\bar{p}=(1/n)\sum_{i=1}^n p_i$.

A theorem of Shepp and Olkin \cite{SO} (see also \cite{M} and \cite{H})
states that
%
\begin{equation}
\label{bi}
H(S)\leq H(\operatorname{bi}(n,\bar{p})),
\end{equation}
where $\operatorname{bi}(n, p)$ denotes the binomial pmf with $n$ trials and
probability $p$ for success. In other words, subject to a fixed mean
$n\bar{p}$, the entropy of $S$ is maximized when all $p_i$ are equal.
Karlin and
Rinott \cite{KR} (see also Harremo\"{e}s \cite{H}) note the
corresponding result
%
\begin{equation}
\label{po}
H(S)\leq H(\operatorname{po}(n\bar{p})),
\end{equation}
where $\operatorname{po}(\lambda)$ denotes the Poisson pmf with mean $\lambda$.

Johnson \cite{J} gives a generalization of (\ref{po}) to \textit{ultra
log-concave} (ULC) distributions. The notion of
ultra log-concavity was introduced by Pemantle \cite{P} in the study of
negative dependence. A pmf $f$
on $\mathbf{Z}_+$ is ULC of order $k$ if $f_i/{{k\choose i}}$ is
log-concave in $i$; it is ULC of order $\infty$,
or simply ULC, if $i!f_i$ is log-concave. Equivalently, these
definitions can be stated with the $\leq_{\mathrm{lc}}$ notation:
\begin{longlist}[1.]
\item[1.]
$f$ is ULC of order $k$ if $f\leq_{\mathrm{lc}} \operatorname{bi}(k, p)$ for some $p\in(0,1)$
(the value of $p$ does not affect the definition);
\item[2.]
$f$ is ULC of order $\infty$ if $f\leq_{\mathrm{lc}} \operatorname{po}(\lambda)$ for some
$\lambda>0$ (the value of $\lambda$ does not affect the definition).
\end{longlist}
An example is the distribution of $S$ in (\ref{po}) and (\ref{bi}).
Denoting the pmf of $S$ by $f^S$, we have
%
\begin{equation}
\label{newton}
f^S \leq_{\mathrm{lc}} \operatorname{bi}(n, \bar{p}),
\end{equation}
which can be shown to be a reformulation of Newton's inequalities (Hardy
\textit{et al.} \cite{HLP}). Also,
note that, as can be verified using the definition, $f$ being ULC of
order $k$
means that it is also ULC of
orders $k+1, k+2, \ldots, \infty$. Another notable property of ULC
distributions,
expressed in our notation, is due to Liggett
\cite{L}.

\begin{theorem}[(\cite{L})]
\label{liggett}
If $f\leq_{\mathrm{lc}} \operatorname{bi}(k,p)$ and $g\leq_{\mathrm{lc}} \operatorname{bi}(m, p), p\in(0,1),$ then
\[
f*g\leq_{\mathrm{lc}} \operatorname{bi}(k+m, p),
\]
where $f*g=\{\sum_{i=0}^j f_i g_{j-i}, j=0, \ldots, k+m\}$ denotes the
convolution of $f$ and $g$.
\end{theorem}

This is a strong result; it implies (\ref{newton}) trivially. Simply observe
that $\operatorname{bi}(1, p_i)\leq_{\mathrm{lc}} \operatorname{bi}(1, \bar{p}),
i=1, \ldots, n,$ and apply
Theorem \ref{liggett}
to obtain $f^S=\operatorname{bi}(1, p_1)*\cdots*\operatorname{bi}(1, p_n)
\leq_{\mathrm{lc}} \operatorname{bi}(n, \bar{p})$,
that is, $f^S$ is ULC of
order $n$. A limiting case of Theorem \ref{liggett} also holds: for
pmf's $f$ and $g$ on
$\mathbf{Z}_+$, we have
\[
f\leq_{\mathrm{lc}} \operatorname{po}(\lambda),\qquad g\leq_{\mathrm{lc}}
 \operatorname{po}  (\mu) \quad\Longrightarrow\quad f*g\leq
_{\mathrm{lc}} \operatorname{po}(\lambda+\mu).
\]

The following generalization of (\ref{po}) is proved by Johnson \cite{J}.
\begin{theorem}[(\cite{J})]
\label{gpo}
If a pmf $f$ on $\mathbf{Z}_+$ is ULC, then
\[
H(f)\leq H(\operatorname{po}(E(f))).
\]
\end{theorem}

Johnson's proof uses two operations, namely convolution with a Poisson
pmf and binomial thinning, to construct a semigroup
action on the set of ULC distributions with a fixed mean. The entropy
is then shown to be monotone along this semigroup. A
corresponding generalization of (\ref{bi}) appears in Yu \cite{Yu}. The
proof adopts the idea of Johnson \cite{J} and is likewise
non-trivial.
\begin{theorem}[(\cite{Yu})]
\label{gbi}
If a pmf $f$ is ULC of order $n$, then
\[
H(f)\leq H\bigl(\operatorname{bi}\bigl(n, E(f)/n\bigr)\bigr).
\]
\end{theorem}

We point out that Theorems \ref{gpo} and \ref{gbi} can be deduced from
Theorem \ref{main}; in fact, both are special cases of the following result.

\begin{theorem}
\label{maxent}
Any log-concave pmf $g$ on $\mathbf{Z}_+$ is the unique maximizer of entropy
in the set $F=\{f\dvt f\leq_{\mathrm{lc}} g,  E(f)=E(g)\}$.
\end{theorem}

\begin{pf}
The log-concavity of $g$ ensures that $\lambda\equiv E(g)<\infty$.
Letting $f\in F$ and using the geometric pmf $\operatorname{ge}(p)=\{p(1-p)^i, i=0,
1,\ldots\}$, we get
\begin{eqnarray*}
D(f|\operatorname{ge}(p))&=&-H(f)-\log(p)-\lambda\log(1-p),
\\
D(g|\operatorname{ge}(p))&=&-H(g)-\log(p)-\lambda\log(1-p),
\end{eqnarray*}
which also shows that $H(f)<\infty$ and $H(g)<\infty$. Since $f\leq
_{\mathrm{lc}} g\leq_{\mathrm{lc}} \operatorname{ge}(p)$, Theorem \ref{main} yields
\[
-H(f)\geq D(f|g)-H(g)\geq-H(g)
\]
so that $H(f)\leq H(g)$ for all $f\in F$, with equality if and only if
$D(f|g)=0$, that is, $f=g$.
\end{pf}

Theorems \ref{gpo} and \ref{gbi} are obtained by noting that both
$\operatorname{po}(\lambda)$ and $\operatorname{bi}(n, p)$ are log-concave. For recent extensions of
Theorems \ref{gpo} and \ref{gbi} to compound distributions, see \cite
{JKM} and \cite{Yu4}.

\subsection{Best binomial approximations in relative entropy}\label{sec2.2}

Recall that $S=\sum_{i=1}^n X_i$ is a sum of independent Bernoulli
random variables, each with success probability $p_i$. Let $\lambda=\sum
_{i=1}^n p_i$ and let $f^S$ denote the pmf of $S$. Approximating $S$
with a Poisson
distribution $\operatorname{Po}(\lambda)$ is an old problem (Le Cam \cite{LC},
Chen \cite{C}, Barbour \textit{et al.} \cite{BHJ}). Approximating $S$
with a binomial $\operatorname{Bi}(n, \bar{p}), \bar{p}=(1/n)\sum_{i=1}^n p_i$,
has also been considered (Stein \cite{S86}, Ehm \cite{E}). The results
are typically
stated in terms of the total variation distance, defined for pmf's $f$
and $g$ as $V(f, g)=\frac{1}{2}\sum_i |f_i-g_i|$.
For example, Ehm \cite{E} applies the method of Stein and Chen to
derive the bound ($\bar{q}=1-\bar{p}$)
\[
V(f^S, \operatorname{bi}(n, \bar{p}))\leq(1-\bar{p}^{n+1}-\bar{q}^{n+1}
)[(n+1)\bar{p}\bar{q}]^{-1} \sum_{i=1}^n
(p_i-\bar{p})^2.
\]

Here, we are concerned with the following problem: what is the best
$m, m\geq n$, and $p\in(0,1)$ for approximating
$S$ with Bi($m, p$)? Intuition says $\operatorname{Bi}(n,\bar{p})$. Indeed, Choi
and Xia \cite{CX} study this in terms of the total
variation distance $d_m=V(f^S, \operatorname{bi}(m, \lambda/m))$ and prove that under
certain conditions, for large enough $m$, $d_m$
increases with $m$.
\begin{theorem}[(\cite{CX})]
\label{choixia}
Let $r=\lfloor\lambda\rfloor$ be the integer part of $\lambda$ and let
$\delta=\lambda-r$. If $r>1+(1+\delta)^2$
and
\[
m\geq\max\bigl\{n, \lambda^2/\bigl(r-1-(1+\delta)^2\bigr)\bigr\},
\]
then $d_m<d_{m+1}<V(f^S, \operatorname{po}(\lambda))$.
\end{theorem}

The derivation of Theorem \ref{choixia} is somewhat involved. However,
if we consider this problem in terms of relative entropy rather than
total variation, then Theorem \ref{triangle} below
gives a definite and equally intuitive answer. Similar results (see
Section~\ref{sec2.3}) hold for the negative
binomial approximation of a sum of independent geometric random variables.
\begin{theorem}
\label{triangle}
Suppose that $m'\geq m\geq n, p'\in(0,1)$. Then,
%
\begin{equation}
\label{tri1}
D(f^S|\operatorname{bi}(m', p'))\geq D\bigl(f^S|\operatorname{bi}(m, \lambda/m)\bigr)
+D\bigl(\operatorname{bi}(m, \lambda/m)|\operatorname{bi}(m', p')\bigr)
\end{equation}
and therefore
\[
D(f^S|\operatorname{bi}(m', p'))\geq D\bigl(f^S|\operatorname{bi}(m, \lambda/m)\bigr)
\geq D(f^S|\operatorname{bi}(n, \bar
{p})).
\]
\end{theorem}

\begin{pf} Let $f=f^S, g=\operatorname{bi}(m, \lambda/m)$ and $h=\operatorname{bi}(m', p')$ in
Theorem \ref{main}. By (\ref{newton}), we
have $f\leq_{\mathrm{lc}} \operatorname{bi}(n,\bar{p})\leq_{\mathrm{lc}} g\leq_{\mathrm{lc}} h$. The claim
follows from (\ref{tri}).
\end{pf}

Theorem \ref{triangle} shows that, for approximating $S$ in the sense
of relative entropy,
\begin{longlist}[1.]
\item[1.]
Bi($m, \lambda/m$), which has the same mean as $S$, is preferable to
Bi($m, p'$), $p'\neq\lambda/m$;
\item[2.]
Bi($n,\bar{p}$) is preferable to Bi($m, \lambda/m$), $m>n$.
\end{longlist}
Obviously, the proof of (\ref{tri1}) still applies when $\operatorname{bi}(m', p')$ is
replaced by $\operatorname{po}(\lambda)$. Hence,
%
\begin{equation}
\label{tripo1}
D(f^S|\operatorname{po}(\lambda))\geq D(f^S|\operatorname{bi}(n,\bar{p}))
+D(\operatorname{bi}(n,\bar{p})|\operatorname{po}(\lambda)),
\end{equation}
that is, Po($\lambda$) is worse than Bi($n,\bar{p})$ by at least
$D(\operatorname{bi}(n,\bar{p})|\operatorname{po}(\lambda))$.

We conclude this subsection with another interesting result in the form
of a corollary of Theorem \ref{main}. Writing $b_m=\operatorname{bi}(m,
\lambda/m)$ for simplicity, we have
\[
D(b_m|\operatorname{po}(\lambda))\geq D(b_m|b_{m+1})+D(b_{m+1}|\operatorname{po}(\lambda))
\]
and, therefore,
%
\begin{equation}
\label{mono}
D(b_m|\operatorname{po}(\lambda))> D(b_{m+1}|\operatorname{po}(\lambda)),\qquad m>\lambda.
\end{equation}
That is, the limit $\operatorname{Bi}(m, \lambda/m)\rightarrow\operatorname{Po}(\lambda
), m\rightarrow\infty,$ is monotone in relative
entropy. As simple as (\ref{mono}) may seem, it is difficult to derive
it directly without Theorem \ref{main}, which
perhaps explains why (\ref{mono}) appears new, even though the
binomial-to-Poisson limit is common knowledge.

\subsection{Analogous results for the negative binomial}\label{sec2.3}

Let $T$ be a sum of geometric random variables, $T=\sum_{i=1}^n Y_i$,
where $Y_i\sim\operatorname{Ge}(r_i)$
independently, $r_i\in(0,1)$. Denote the mean of $T$ by $\mu=\sum
_{i=1}^n (1-r_i)/r_i$ and denote the pmf of
$T$ by $f^T$. Let $\operatorname{nb}(n, r)=\{{{n+i-1}\choose{i}} r^n(1-r)^i, i=0,
1,\ldots\}$ denote the pmf of the negative binomial $\operatorname{NB}(n,r$).

The counterpart of (\ref{bi}) appears in Karlin and Rinott \cite{KR}.
\begin{theorem}[(\cite{KR})]
\label{nb} $H(T)\geq H(\operatorname{nb}(n, n/(n+\mu))).$
\end{theorem}

In other words, subject to a fixed mean $\mu$, the entropy of $T$ is
minimized when all $r_i$ are equal. Theorem \ref{nb}
can be generalized as follows.
\begin{theorem}
\label{nmaxent}
Any log-concave pmf $f$ is the unique minimizer of entropy in the set
$G=\{g\dvt f\leq_{\mathrm{lc}} g\leq_{\mathrm{lc}} \operatorname{ge}(p),
E(g)=E(f)\}, p\in(0,1)$.
\end{theorem}

We realize that Theorem \ref{nmaxent} is just a reformulation of
Theorem~\ref{maxent}, which follows from Theorem~\ref{main}. To show
that Theorem
\ref{nmaxent} indeed implies Theorem \ref{nb}, we need the following
inequality of Hardy \textit{et al.}
\cite{HLP}, written in our notation as
%
\begin{equation}
\label{hlp}
\operatorname{nb}\bigl(n, n/(n+\mu)\bigr)\leq_{\mathrm{lc}} f^T.
\end{equation}
We also need $f^T$ to be log-concave, but this holds because
convolutions of log-concave sequences are also log-concave.

Next, we consider the problem of selecting the best $m, m\geq n$, and
$r\in(0,1)$ for
approximating $T$ with NB($m, r$).

\begin{theorem}
\label{TRI2}
Suppose $m'\geq m\geq n$ and $r'\in(0,1)$. Write $\operatorname{nb}_m=\operatorname{nb}(m, m/(m+\mu
))$ as shorthand. Then,
\[
D(f^T|\operatorname{nb}(m', r'))\geq D(f^T|\operatorname{nb}_m)+D(\operatorname{nb}_m|\operatorname{nb}(m', r'))
\]
and, therefore,
\[
D(f^T|\operatorname{nb}(m', r'))\geq D(f^T|\operatorname{nb}_m)\geq D\bigl(f^T|\operatorname{nb}\bigl(n, n/(n+\mu)\bigr)\bigr).
\]
\end{theorem}

\begin{pf} The relations
\[
\operatorname{nb}(m', r')\leq_{\mathrm{lc}} \operatorname{nb}_m
\leq_{\mathrm{lc}} \operatorname{nb}\bigl(n, n/(n+\mu)\bigr)
\]
are easy to verify. We also have (\ref{hlp}). The claim follows from
(\ref{tri-}).
\end{pf}

Theorem \ref{TRI2} implies that for approximating $T$ in the sense of
relative entropy,
NB($n,n/(n+\mu)$) is no worse than NB($m', r'$) whenever $m'\geq n$.
The counterpart of (\ref{tripo1}) also holds
($\operatorname{nb}_n= \operatorname{nb}(n, n/(n+\mu))$):
\[
D(f^T|\operatorname{po}(\mu))\geq D(f^T|\operatorname{nb}_n)+D(\operatorname{nb}_n|\operatorname{po}(\mu)),
\]
that is, Po($\mu$) is worse than NB($n,n/(n+\mu)$) by at least
$D(\operatorname{nb}_n|\operatorname{po}(\mu))$.

In addition, parallel to (\ref{mono}), we have
%
\begin{equation}
\label{nmono}
D(\operatorname{nb}_m|\operatorname{po}(\mu))
> D(\operatorname{nb}_{m'}|\operatorname{po}(\mu)),\qquad m'>m>0,
\end{equation}
that is, the limit $\operatorname{NB}(m, m/(m+\mu))\rightarrow\operatorname{Po}(\mu),\
m\rightarrow\infty,$ is monotone in relative entropy. Note
that in (\ref{nmono}), $m$ and $m'$ need not be integers; similarly in
Theorem \ref{TRI2}.

We conclude this subsection with a problem on the behavior of $\leq
_{\mathrm{lc}}$ under convolution. Analogous to
Theorem \ref{liggett} is the following result of Davenport and P\'
{o}lya (\cite{DP}, Theorem 2), rephrased in terms of
$\leq_{\mathrm{lc}}$.
\begin{theorem}[(\cite{DP})]
\label{conj}
Suppose that pmf's $f$ and $g$ on $\mathbf{Z}_+$ satisfy $\operatorname{nb}(k, r)\leq
_{\mathrm{lc}} f, \operatorname{nb}(m, r)\leq_{\mathrm{lc}} g$ for $k, m>0, r\in(0, 1)$. Their
convolution $f*g$ then satisfies
\[
\operatorname{nb}(k+m, r)\leq_{\mathrm{lc}} f*g.
\]
\end{theorem}

Actually, Davenport and P\'{o}lya \cite{DP} assume that $k+m=1$, so
their conclusion is the log-convexity of $f*g$, but it is readily
verified that the same proof works for all positive $k$ and $m$. The
limiting case also holds, that is,
\[
\operatorname{po}(\lambda)\leq_{\mathrm{lc}} f,\qquad
 \operatorname{po}(\mu)\leq_{\mathrm{lc}} g \quad \Longrightarrow\quad \operatorname{po}(\lambda
+\mu)\leq_{\mathrm{lc}} f*g.
\]
An open problem is to determine general conditions that ensure
%
\begin{equation}
\label{open}
f\leq_{\mathrm{lc}} f',\qquad g\leq_{\mathrm{lc}} g'\quad \Longrightarrow\quad f*g\leq_{\mathrm{lc}} f'*g'.
\end{equation}
Theorem \ref{liggett} simply says that (\ref{open}) holds if
$f'=\operatorname{bi}(k,p)$ and $g'=\operatorname{bi}(m, p)$ with the same $p$ and Theorem \ref{conj}
says that (\ref{open}) holds if $f=\operatorname{nb}(k, r)$ and $g=\operatorname{nb}(m, r)$ with the
same $r$. The proofs of Theorems \ref{conj} and \ref{liggett} (Theorem
\ref{liggett} especially) are non-trivial. It is reasonable to ask
whether there exist other interesting and non-trivial instances of (\ref{open}).

\section[Proof of Theorem 1]{Proof of Theorem \protect\ref{main}}\label{sec3}

The proof of Theorem \ref{main} hinges on the following lemma that
dates back to Karlin and Novikoff \cite{KN} and Karlin
and Studden \cite{KS}. Our assumptions are slightly different from
those of Karlin and Studden \cite{KS}, Lemma XI. 7.2. In the proof
(included for completeness), the number of sign changes of a sequence
is counted discarding zero terms.

\begin{lemma}[(\cite{KS})]
\label{pmp}
Let $a_i, i=0, 1, \ldots,$ be a real sequence such that $\sum
_{i=0}^\infty a_i=0$ and $\sum_{i=0}^\infty i\times a_i=0.$
Suppose that the set $C=\{i\dvt a_i>0\}$ is an interval on $\mathbf
{Z}_+$. For any concave function
$w(i)$ on $\mathbf{Z}_+$, we then have
%
\begin{equation}
\label{wcv}
\sum_{i=0}^\infty w(i)a_i\geq0.
\end{equation}
\end{lemma}

\begin{pf} Karlin and Studden (\cite{KS}, Lemma XI. 7.2) assume that
$a_i, i=0, 1, \ldots,$ changes sign exactly twice, with sign
sequence $-, +, -$. However, it also suffices to assume that $C$ is an
interval. Suppose that $a_i$ changes sign exactly
once, with sign sequence $+, -$, that is, there exists $0\leq k<\infty$
such that $a_i\geq0, 0\leq i\leq k$, with strict
inequality for at least one $i\leq k$, and $a_i\leq0, i>k$. Then,
\[
\sum_{i=0}^\infty ia_i \leq\sum_{i=0}^k ka_i +\sum_{i=k+1}^\infty
(k+1) a_i =-\sum_{i=0}^k a_i<0,
\]
a contradiction. Similarly, the sign sequence cannot be $-, +$ either.
Assuming that $C$ is an interval,
this shows that, except for the trivial case $a_i\equiv0$, the
sequence $a_i$ changes sign exactly twice, with sign sequence
$-, +, -$.

The rest of the argument is well known. We proceed to show that the
sequence $A_j=\sum_{i=0}^j a_i$ has exactly one sign
change, with sign sequence $-, +$. Similarly, $\sum_{i=0}^j A_i\leq0$
for all $j=0, 1, \ldots,$ which implies (\ref{wcv}) for
every concave function $w(i)$ upon applying summation by parts.
\end{pf}

Theorem \ref{concave} below is a consequence of Lemma \ref{pmp}.
Although not phrased as such, the basic idea is
implicit in Karlin and Studden \cite{KS} in their analyses of special
cases; see also Whitt \cite{W}. When $f$ is the pmf of a
sum of $n$ independent Bernoulli random variables and $g=\operatorname{bi}(n,
E(f)/n),$ as discussed in Section 2, Theorem \ref{concave} reduces to an
inequality of Hoeffding \cite{H56}.
\begin{theorem}
\label{concave}
Suppose that two pmf's $f$ and $g$ on $\mathbf{Z}_+$ satisfy $f\leq
_{\mathrm{lc}} g$ and $E(f)=E(g)<\infty.$ For any concave function $w(i)$ on
$\mathbf{Z}_+$, we then have
\[
\sum_{i=0}^\infty f_iw(i)\geq\sum_{i=0}^\infty g_i w(i).
\]
\end{theorem}

\begin{pf} Since $E(g)<\infty$ and $w$ is concave, $\sum_{i=0}^\infty
g_i w(i)$ either converges absolutely or diverges to
$-\infty$. Assume the former. Since $\log(f_i/g_i)$ is concave and
hence unimodal, the set $C=\{i\dvt f_i-g_i>0\}$ must be an
interval. The result then follows from Lemma \ref{pmp}.
\end{pf}

Theorem \ref{main} is a consequence of Theorem \ref{concave}. Actually,
we prove a slightly more general
``quadrangle inequality,'' which may be of interest. Theorem \ref{main}
corresponds to the special case $g=g'$ in Theorem \ref{main4}.

\begin{theorem}
\label{main4}
Let $f, g, g', h$ be pmf's on $\mathbf{Z}_+$ such that $f\leq_{\mathrm{lc}} g\leq
_{\mathrm{lc}} g'\leq_{\mathrm{lc}} h$. If
$E(f)=E(g)<\infty$, then $D(f|h)<\infty$ and
%
\begin{equation}
\label{tri4}
D(f|h)+D(g|g')\geq D(f|g')+D(g|h);
\end{equation}
if $E(g')=E(h)<\infty$, then
%
\begin{equation}
\label{tri4-}
D(h|f)+D(g'|g)\geq D(g'|f)+D(h|g).
\end{equation}
\end{theorem}

\begin{pf} The concavity of $\log(f_i/h_i)$ and $E(f)<\infty$ imply
$D(f|h)<\infty$. Likewise for $D(g|h)$. Thus, (\ref{tri4}) can be
written as
\[
D(f|h)-D(f|g')\geq D(g|h)-D(g|g')
\]
or, equivalently,
%
\begin{equation}
\label{ineq4}
\sum_i f_i \log(g'_i/h_i)\geq\sum_i g_i\log(g'_i/h_i).
\end{equation}
Since $\log(g'_i/h_i)$ is concave in $\operatorname{supp}(g')$, and $\operatorname{supp}(f)\subset
\operatorname{supp}(g)\subset \operatorname{supp}(g')$, (\ref{ineq4}) follows directly
from Theorem \ref{concave}.

To prove (\ref{tri4-}), we may assume $D(h|f)<\infty$ and
$D(g'|g)<\infty$. These imply, in particular, that
$\operatorname{supp}(f)=\operatorname{supp}(g')=\operatorname{supp}(h)$. We get
\[
\sum_i g'_i\log(f_i/g_i)\geq\sum_i h_i\log(f_i/g_i)
\]
and (\ref{tri4-}) follows as before.
\end{pf}

\section{The continuous case}\label{sec4}

For probability density functions (pdf's) $f$ and $g$ with respect to
Lebesgue measure on $\mathbf{R}$, the differential entropy of $f$ and
the relative entropy between $f$ and $g$ are defined, respectively, as
\[
H(f)=\int_{-\infty}^\infty-f(x)\log(f(x)) \, \mathrm{d}x \quad\mbox{and}\quad
D(f|g)=\int_{-\infty}^\infty f(x)\log\bigl(f(x)/g(x)\bigr) \, \mathrm{d}x.
\]
Parallel to the discrete case, let us write $f\leq_{\mathrm{lc}} g$ if
\begin{longlist}
\item[1.]
$\operatorname{supp}(f)$ and $\operatorname{supp}(g)$ are both intervals on $\mathbf{R}$;
\item[2.]
$\operatorname{supp}(f)\subset \operatorname{supp}(g)$; and
\item[3.]
$\log(f(x)/g(x))$ is concave in $\operatorname{supp}(f)$.
\end{longlist}
There then holds a continuous analog of Theorem \ref{main} (with its
first phrase replaced by ``Let $f, g, h$ be pdf's on $\mathbf{R}$'');
the proof is similar and is hence omitted.

The following maximum/minimum entropy result parallels Theorems \ref
{maxent} and \ref{nmaxent}.

\begin{theorem}
\label{minicont}
If a pdf $g$ on $\mathbf{R}$ is log-concave, then it maximizes the
differential entropy in the set $F=\{f\dvt f\leq_{\mathrm{lc}} g,
E(f)=E(g)\}$. Alternatively, if a pdf $f$ on $\mathbf{R}$ is
log-concave, then it minimizes the differential entropy in the
set $G=\{g\dvt f\leq_{\mathrm{lc}} g, g \mbox{ is log-concave and } E(g)=E(f)\}$.
\end{theorem}

We illustrate Theorem \ref{minicont} with a minimum entropy
characterization of the gamma distribution. This parallels
Theorem \ref{nb} for the negative binomial. Denote by $\operatorname{gam}(\alpha, \beta
)$ the pdf of the gamma distribution Gam$(\alpha, \beta)$, that is,
\[
\operatorname{gam}(x; \alpha, \beta)=\beta^{-\alpha} x^{\alpha-1} \mathrm{e}^{-x/\beta}/\Gamma
(\alpha),\qquad x>0.
\]

\begin{theorem}
\label{gam}
Let $\alpha_i\geq1, \beta_i>0$ and let $X_i\sim \operatorname{Gam}(\alpha_i,  1),\
i=1, \ldots, n,$ independently. Define
$S=\sum_{i=1}^n \beta_i X_i$. Then, subject to a fixed mean $ES=\sum
_{i=1}^n \alpha_i \beta_i$, the differential entropy of $S$ $($as
a function of $\beta_i, i=1,\ldots, n)$ is minimized when all $\beta
_i$ are equal.
\end{theorem}

Note that Theorem 3.1 of Karlin and Rinott (\cite{KR}; see also Yu \cite
{Yu3}) implies that Theorem
\ref{gam} holds when all $\alpha_i$ are equal. We use $\leq_{\mathrm{lc}}$ to
give an extension to general $\alpha_i\geq1$. A useful result is Lemma
\ref{gamconj}, which reformulates Theorem 4 of Davenport and P\'{o}lya
\cite{DP}. As in Theorem \ref{conj}, Davenport and P\'{o}lya assume
$\alpha_1+\alpha_2=1$, but the proof works for all positive $\alpha_1,
\alpha_2$.

\begin{lemma}[(\cite{DP}, Theorem 4)]
\label{gamconj}
Let $\alpha_1, \alpha_2>0$ and let $f$ and $g$ be pdf's on $(0, \infty
)$ such that $\operatorname{gam}(\alpha_1, 1)\leq_{\mathrm{lc}} f$ and $\operatorname{gam}(\alpha_2, 1)\leq
_{\mathrm{lc}} g$. Then,
\[
\operatorname{gam}(\alpha_1+\alpha_2, 1)\leq_{\mathrm{lc}} f*g,
\]
where $(f*g)(x)=\int_0^x f(y)g(x-y) \, \mathrm{d}y$.
\end{lemma}

\begin{pf*}{Proof of Theorem \protect\ref{gam}} Repeated application of Lemma \ref
{gamconj} yields
%
\begin{equation}
\label{gam0}
\operatorname{gam}(\alpha_+, 1)\leq_{\mathrm{lc}} f^S,
\end{equation}
where $\alpha_+=\sum_{i=1}^n \alpha_i$ and $f^S$ denotes the pdf of
$S$. Alternatively, we can show (\ref{gam0}) by noting that $f^S$ is a
mixture of $\operatorname{gam}(\alpha_+, \beta)$, where $\beta$ has the distribution
of $S/\sum_{i=1}^n X_i$ (see,\vspace{1pt} e.g., \cite{W} and \cite{Yu2}).
Since $\alpha_i\geq1$, each $X_i$ is log-concave and so is $f^S$. The
claim follows from Theorem~\ref{minicont}.
\end{pf*}

Weighted sums of gamma variates, as in Theorem \ref{gam}, arise
naturally in statistical contexts, for example, as quadratic forms in
normal variables, but their distributions can be non-trivial to compute
(Imhof \cite{I}). When comparing different gamma distributions as
convenient approximations, we obtain a result similar to Theorems \ref
{triangle} and \ref{TRI2}. The proof, also similar, is omitted.

\begin{theorem}
\label{gam2}
Fix $\alpha_i>0, \beta_i>0$ and let $X_i\sim \operatorname{Gam}(\alpha_i, 1), i=1,
\ldots, n,$ independently. Define $S=\sum_{i=1}^n \beta_i X_i$, with
pdf $f^S$. Write $g_a=\operatorname{gam}(a, \sum_{i=1}^n \beta_i\alpha_i/a)$ as
shorthand. For $b>0$ and $a'\geq a\geq\alpha_+$, where $\alpha_+=\sum
_{i=1}^n \alpha_i$, we then have
\[
D(f^S|\operatorname{gam}(a', b))\geq D(f^S|g_a)+D(g_a|\operatorname{gam}(a', b))
\]
and, consequently,
\[
D(f^S|\operatorname{gam}(a', b))\geq D(f^S|g_a)\geq D(f^S|g_{\alpha_+}).
\]
\end{theorem}

In other words, to approximate $S$ in the sense of relative entropy,
Gam$(\alpha_+, \sum_{i=1}^n \beta_i\alpha_i/\alpha_+)$,
which has the same mean as $S$, is no worse than Gam$(a, b)$ whenever
$a\geq\alpha_+$. Note that, unlike in Theorem
\ref{gam}, we do not require here that $\alpha_i\geq1$.

Overall, there is a remarkable parallel between the continuous and
discrete cases.

\section*{Acknowledgments}
The author would like to thank three referees for their constructive comments.

\printhistory

\end{document}